\documentclass[12pt]{amsart}

\usepackage{fancyhdr}

\usepackage{amssymb,latexsym,amscd}

\usepackage{amsthm,amsmath}

\input xy

\xyoption{all}

\usepackage[cp850]{inputenc}

\newcommand{\cO}{\mathcal{O}}

\newcommand{\cE}{\mathcal{E}}

\newcommand{\cM}{\mathcal{M}}

\newcommand{\cB}{\mathcal{B}}

\newcommand{\cL}{\mathcal{L}}

\newcommand{\Sym}{\mathrm{Sym}}
\newcommand{\Pic}{\mathrm{Pic}}

\newcommand{\lra}{\longrightarrow}

\newcommand{\PP}{\mathbb{P}}

\newcommand{\Ext}{\mathrm{Ext}}

\newcommand{\im}{\mathrm{im}}
\newcommand{\rk}{\mathrm{rk}}

\def\map#1{\ \smash{\mathop{\longrightarrow}\limits^{#1}}\ }

\theoremstyle{plain}

\newtheorem{thm}{Theorem}[section]

\newtheorem{prop}[thm]{Proposition}

\begin{document}

\title[]{On the base locus of the linear system of generalized theta functions}

\author{Christian Pauly}

\address{D\'epartement de Math\'ematiques \\ Universit\'e de Montpellier II - Case Courrier 051 \\ Place Eug\`ene Bataillon \\ 34095 Montpellier Cedex 5 \\ France}
\email{pauly@math.univ-montp2.fr}




\begin{abstract}
Let $\cM_r$ denote the moduli space of semi-stable rank-$r$ vector bundles 
with trivial determinant over a smooth projective curve $C$ of genus $g$.
In this paper we study the base locus $\cB_r \subset \cM_r$ of the linear
system of the determinant line bundle $\cL$ over $\cM_r$, i.e., the set 
of semi-stable rank-$r$ vector bundles without theta divisor. We construct base points in 
$\cB_{g+2}$ over any curve $C$, and base points in $\cB_4$ over any
hyperelliptic curve. 
\end{abstract}

\maketitle


\section{Introduction}
Let $C$ be a complex smooth projective curve of genus $g$ and let $\cM_r$  denote the
coarse moduli space parametrizing semi-stable rank-$r$ vector bundles with trivial 
determinant over the curve $C$. Let $\cL$ be the determinant line bundle over the
moduli space $\cM_r$ and let $\Theta \subset \Pic^{g-1}(C)$ be the Riemann theta
divisor in the degree $g-1$ component of the Picard variety of $C$. By \cite{BNR}
there is a canonical isomorphism 
$|\cL|^* \map{\sim} |r\Theta|$, under which the natural rational map 
$\varphi_\cL : \cM_r \dasharrow |\cL|^*$ is identified with the so-called theta map
$$\theta : \cM_r \dasharrow |r \Theta|, \qquad E \mapsto \theta(E) \subset 
\Pic^{g-1}(C).$$
The underlying set of $\theta(E)$ consists of line bundles $L \in \Pic^{g-1}(C)$
with $h^0(C, E \otimes L ) > 0$. For a general semi-stable vector bundle $E$, $\theta(E)$
is a divisor. If $\theta(E) = \Pic^{g-1}(C)$, we say that $E$ has no theta 
divisor. We note that the indeterminacy locus
of the theta map $\theta$, i.e., the set of bundles $E$ without theta divisor, 
coincides with the base locus $\cB_r \subset \cM_r$ of the linear system $|\cL|$. 

\bigskip

Over the past years many authors \cite{A}, \cite{B2},  \cite{He}, \cite{Hi}, \cite{P}, \cite{R}, \cite{S} have studied the base locus $\cB_r$ of $|\cL|$ and their analogues for the 
powers $|\cL^k|$. For a recent survey of this subject we refer to \cite{B1}.

\bigskip

It is natural to introduce for a curve $C$ the integer $r(C)$ defined as 
the minimal rank for which there exists a semi-stable rank-$r(C)$ 
vector bundle with trivial determinant over $C$ without theta divisor (see also \cite{B1}
section 6). It is known \cite{R} that $r(C) \geq 3$ for any curve $C$ and 
that $r(C) \geq 4$ for a
generic curve $C$. Our main result shows the existence of vector bundles of low ranks
without theta divisor.

\begin{thm}
We assume that $g \geq 2$. Then we have the following bounds.
\begin{enumerate}
\item $r(C) \leq g+2$.
\item $r(C) \leq 4$, if $C$ is hyperelliptic. 
\end{enumerate}
\end{thm}

The first part of the theorem improves the upper bound $r(C) \leq \frac{(g+1)(g+2)}{2}$
given in \cite{A}. The statements of the theorem are equivalent to the existence of a 
semi-stable rank-$(g+2)$ (resp. rank-$4$) vector bundle without theta divisor --- see
section 2.1 (resp. 2.2). The construction of these vector bundles uses ingredients which are already implicit in \cite{Hi}.

\bigskip

Theorem 1.1 seems to hint towards a dependence of the integer $r(C)$ on the
curve $C$.

\bigskip

{\em Notations:} If $E$ is a vector bundle over $C$, we will write $H^i(E)$ for $H^i(C,E)$
and $h^i(E)$ for $\dim H^i(C,E)$. We denote the slope of $E$ by $\mu(E) := 
\frac{\deg E}{\rk E}$, the canonical bundle over $C$ by $K$ and the degree $d$
component of the Picard variety of $C$ by $\Pic^d(C)$. 

\bigskip

\section{Proof of Theorem 1.1}

\subsection{Semi-stable rank-$(g+2)$ vector bundles without theta divisor}

We consider a line bundle $L \in \Pic^{2g+1}(C)$. Then $L$ is globally generated,
$h^0(L) = g+2$ and the evaluation bundle $E_L$, which is defined by the exact
sequence
\begin{equation} \label{es1}
0 \lra E^*_L \lra H^0(L) \otimes \cO_C \map{ev} L \lra 0,
\end{equation}
is stable (see e.g. \cite{Bu}), with $\deg E_L = 2g+1$, $\rk E_L = g+1$ and $\mu(E_L) = 2 -
\frac{1}{g+1}$.

\bigskip

A cohomology class $e \in \Ext^1(LK^{-1}, E_L) = H^1(E_L \otimes KL^{-1})$ determines
a rank-$(g+2)$ vector bundle $\cE_e$ given as an extension
\begin{equation} \label{es2}
0 \lra E_L \lra \cE_e \lra LK^{-1} \lra 0.
\end{equation}

\begin{prop}
For any non-zero class $e$, the rank-$(g+2)$ vector bundle $\cE_e$ is semi-stable. 
\end{prop}

\begin{proof}
Consider a proper subbundle $A \subset \cE_e$. If $A \subset E_L$, then $\mu(A) \leq
\mu(E_L) = 2 - \frac{1}{g+1}$ by stability of $E_L$, so the subbundles of $E_L$ 
cannot destabilize $\cE_e$. If $A \nsubseteq E_L$, we introduce 
$S = A \cap E_L \subset E_L$ and consider the exact sequence 
$$0 \lra S \lra A \lra LK^{-1}(-D) \lra 0,$$
where $D$ is an effective divisor. If $\rk S = g+1$ or $S = 0$, we easily conclude 
that $\mu(A) < \mu(\cE_e) =2$. If $\rk S < g+1$ and $S \not= 0$, then stability of $E_L$ gives the inequality
$\mu(S) < \mu(E_L) = 2 - \frac{1}{g+1}$. We introduce the integer 
$\delta  = 2 \rk S - \deg S$. Then the previous inequality is equivalent to $\delta \geq 1$.
Now we compute
$$ \mu(A) = \frac{\deg S + \deg LK^{-1}(-D)}{\rk S +1} \leq \frac{2 \rk S - \delta + 3}{\rk S +1}
= 2 + \frac{1-\delta}{\rk S +1} \leq 2 = \mu(\cE_e),$$
which shows the semi-stablity of $\cE_e$.
\end{proof}

We tensorize the exact sequence \eqref{es1} with $L$ and take the cohomology

\begin{equation} \label{esEL}
0 \lra H^0(E_L^* \otimes L) \lra H^0(L) \otimes H^0(L) \map{\mu} H^0(L^2) \lra 0.
\end{equation}

Note that $h^1(E_L^* \otimes L ) = h^0(E_L \otimes KL^{-1}) = 0$ by stability of $E_L$. The
second map $\mu$ is the multiplication map and factorizes through $\Sym^2 H^0(L)$, i.e.,
$$ \Lambda^2 H^0(L) \subset H^0(E_L^* \otimes L) = \ker \mu.$$
By Serre duality a cohomology class $e \in \Ext^1(LK^{-1}, E_L) = H^1(E_L \otimes KL^{-1}) =
H^0(E_L^* \otimes L)^*$ can be viewed as a hyperplane $H_e \subset H^0(E_L^* \otimes L)$. Then we have the following

\begin{prop}
If $\Lambda^2 H^0(L) \subset H_e$, then the vector bundle $\cE_e$ satisfies 
$$ h^0(\cE_e \otimes \lambda) > 0, \qquad \forall \lambda \in \Pic^{g-3}(C).$$
\end{prop}

\begin{proof}
We tensorize the exact sequence \eqref{es2} with $\lambda \in \Pic^{g-3}(C)$ and
take the cohomology
$$ 0 \lra H^0(E_L \otimes \lambda) \lra H^0(\cE_e \otimes \lambda) \lra
H^0(LK^{-1} \lambda) \map{\cup e} H^1(E_L \otimes \lambda) \lra \cdots $$
Since $\deg LK^{-1} \lambda = g$, we can write $LK^{-1} \lambda = \cO_C(D)$ for 
some effective divisor $D$. It is enough to show that $h^0(\cE_e \otimes \lambda) > 0$
holds for $\lambda$ general. Hence we can assume that $h^0(LK^{-1}\lambda) = 
h^0(\cO_C(D)) = 1$.

\bigskip

If $h^0(E_L \otimes \lambda) >0$, we are done. So we assume 
$h^0(E_L \otimes \lambda) = 0$, which implies $h^1(E_L \otimes \lambda) =1$ by 
Riemann-Roch. Hence we obtain that $h^0(\cE_e \otimes \lambda) > 0$ if and only if
the cup product map 
$$ \cup e : H^0(\cO_X(D)) \lra H^1(E_L \otimes \lambda) = 
H^0(E^*_L \otimes L(-D))^* $$
is zero. Furthermore $\cup e$ is zero if and only if $H^0(E^*_L \otimes L(-D))
\subset H_e$. Now we will show the inclusion 
\begin{equation} \label{inclusionD}
H^0(E^*_L \otimes L(-D)) \subset \Lambda^2 H^0(L).
\end{equation}
We tensorize the exact sequence \eqref{es1} with $L(-D)$ 
and take cohomology 
$$ 0 \lra H^0(E_L^* \otimes L(-D)) \lra H^0(L) \otimes H^0(L(-D)) \map{\mu} H^0(L^2(-D)) 
\lra \cdots$$
Since $h^0(E_L^* \otimes L(-D)) = 1$, we conclude that $h^0(L(-D)) = 2$ and 
$H^0(E^*_L \otimes L(-D)) = \Lambda^2 H^0(L(-D)) \subset \Lambda^2 H^0(L)$.

\bigskip
Finally the proposition follows: if $\Lambda^2 H^0(L) \subset H_e$, then by \eqref{inclusionD}
$H^0(E^*_L \otimes L(-D)) \subset H_e$ for general $D$, or equivalently 
$h^0(\cE_e \otimes \lambda) > 0 $ for general $\lambda \in \Pic^{g-3}(C)$.
\end{proof}

We introduce the linear subspace $\Gamma \subset \Ext^1(LK^{-1}, E_L)$ defined by
$$ \Gamma := \ker \left(
\Ext^1(LK^{-1}, E_L) = H^0(E_L^* \otimes L)^* \lra 
\Lambda^2 H^0(L)^* \right),$$
which has dimension $\frac{g(g-1)}{2} > 0$. Then for any non-zero
cohomology class $e \in \Gamma$ and any $\gamma \in \Pic^2(C)$ satisfying 
$\gamma^{g+2} = L^2 K^{-1} = \det \cE_e$, the rank-$(g+2)$ vector bundle  
$$ \cE_e \otimes \gamma^{-1} $$
has trivial determinant, is semi-stable by Proposition 2.1 and has no theta
divisor by Proposition 2.2.

\subsection{Hyperelliptic curves}

In this subsection we assume that $C$ is hyperelliptic and we denote by $\sigma$ the hyperelliptic
involution. The construction of a semi-stable rank-$4$ vector bundle without theta
divisor has been given in \cite{Hi} section 6 in the case $g=2$, but it can be
carried out for any $g \geq 2$ without major modification. For the convenience of the
reader, we recall the construction and refer to \cite{Hi} for the details and
the proofs.

\bigskip

Let $w \in C$ be a Weierstrass point. Any non-trivial extension 
$$ 0 \lra \cO_C(-w) \lra G \lra \cO_C \lra 0$$
is a stable, $\sigma$-invariant, rank-$2$ vector bundle with $\deg G = -1$. 
By \cite{Hi} Theorem 4 a cohomology class $e \in H^1(\Sym^2 G)$ determines 
a symplectic rank-$4$ bundle
$$ 0 \lra G \lra \cE_e \lra G^* \lra 0.$$
Moreover it is easily seen that, for any non-zero class $e$, the vector bundle
$\cE_e$ is semi-stable. By \cite{Hi} Lemma 16 the composite map
$$ D_G : \PP H^1 (\Sym^2 G) \lra \cM_4 \map{\theta} |4 \Theta|, \qquad 
e \mapsto \theta(\cE_e)$$
is the projectivization of a linear map
$$ \widetilde{D_G} :  H^1 (\Sym^2 G) \lra H^0(\Pic^{g-1}(C), 4 \Theta).$$
The involution $i(L) = KL^{-1}$ on $\Pic^{g-1}(C)$ induces a linear involution on
$|4\Theta|$ with eigenspaces $|4\Theta|_\pm$. Note that $4\Theta \in |4\Theta|_+$.
We now observe that $\theta(\cE) \in |4\Theta|_+$ for any symplectic rank-$4$ vector
bundle $\cE$ --- see e.g. \cite{B2}. Moreover we have the equality 
$\theta(\sigma^* \cE) = i^* \theta(\cE)$ for any vector bundle $\cE$. These two
observations imply that the linear map $\widetilde{D_G}$ is equivariant with 
respect to the induced involutions $\sigma$ and $i$. Since $\im  \ \widetilde{D_G}
\subset H^0(\Pic^{g-1}(C), 4 \Theta)_+$, we obtain that one of the two eigenspaces
$H^1(\Sym^2 G)_\pm$ is contained in the kernel $\ker \widetilde{D_G}$, hence
give base points for the theta map. We now compute as in \cite{Hi} using the
Atiyah-Bott-fixed-point formula
$$ h^1(\Sym^2 G)_+ = g-1, \qquad h^1(\Sym^2 G)_- = 2g+1.$$
One can work out that $H^1(\Sym^2 G)_+ \subset \ker \widetilde{D_G}$. Hence any
$\cE_e$ with non-zero $e \in H^1(\Sym^2 G)_+$ is a semi-stable rank-$4$ vector bundle without
theta divisor.

\end{document}